\documentclass{amsart}
\usepackage{amssymb, mathrsfs, bbm}
\usepackage[sc, osf]{mathpazo}
\selectfont
\usepackage[protrusion=true, expansion=true]{microtype}

\usepackage[pdftex]{hyperref}

\swapnumbers
\theoremstyle{plain}

\newtheorem{scholium}[subsection]{Scholium}

\theoremstyle{remark}

\theoremstyle{definition}

\numberwithin{equation}{subsection}

\def\11{\mathbf{1}}

\def\AA{\mathbf{A}}

\def\FF{\mathbf{F}} 
\def\GG{\mathbf{G}}

\def\ZZ{\mathbf{Z}}

 \def\dim{\mathrm{dim}}

 \def\gr{\mathrm{gr}}
 \def\Gr{\mathcal{G}r}

\makeatletter
\renewcommand{\@makefnmark}{\mbox{\textsuperscript{}}}
\makeatother

\title{On Hall-Littlewood polynomials}
\author{R. Virk}
\email{rsvirk@gmail.com}
\begin{document}
\maketitle
\subsection{Introduction}
This billet should be regarded as a footnote to \cite{GL}. We observe that Hall-Littlewood polynomials encode motivic Euler characteristics of intersections in the affine Grassmanian. Consequently, some coarse Hodge theoretic information is contained in them.

The reader should view this as an exercise in semantics. After all, Hall-Littlewood polynomials have traditionally been interpreted in terms of point counting (over finite fields). The Hodge theoretic interpretation (as opposed to point counting) is merely a preference for the Hodge realization (as opposed to the $\ell$-adic realization) of motivic cohomology.
Of course, all of these cohomological interpretations are, in turn, a rephrasing of I.G. Macdonald's fantastic investigations of Hall-Littlewood polynomials in connection with $p$-adic groups (the spherical functions of \cite{Mac}).

Embarrassingly, your correspondent has arrived at this understanding at a rather glacial pace. Hopefully, writing some of this down will speed things up for others encountering similar esoterica.

\subsection{The setup}
Work over the complex numbers. Fix a connected reductive linear algebraic group, a Borel subgroup, and a maximal torus. To this data one associates:
\begin{align*}
X &- \mbox{the group of coweights}, \\
X_+ &- \mbox{the semi-group of dominant coweights}, \\
\langle -, - \rangle &- \mbox{the evident perfect pairing between weights and coweights}, \\
\Gr &- \mbox{the affine Grassmanian}.
\end{align*}
The affine Grassmanian is an ind-projective ind-variety. Details on the construction/properties of $\Gr$ and subvarieties appearing below can be found in \cite{PS} and \cite{MV}. The notation used in this note is chosen to conform with \cite{MV}.

The disk group associated to our reductive group acts on $\Gr$. The orbits of this action are parameterized by $X_+$. Write $\Gr^{\lambda}$ for the orbit corresponding to $\lambda\in X_+$. The $\Gr^{\lambda}$ are locally closed (finite dimensional) subvarieties of $\Gr$. The affine Grassmanian also decomposes into orbits for the loop group of the unipotent radical of the Borel subgroup opposite to our chosen Borel. These orbits are parameterized by $X$. Write $S_{\nu}$ for the orbit corresponding to $\nu\in X$.
\begin{align*}
\Gr &= \bigsqcup_{\lambda\in X_+}\Gr^{\lambda} && \mbox{(Cartan decomposition)}, \\
\Gr &= \bigsqcup_{\nu\in X}S_{\nu} && \mbox{(Iwasawa decomposition)}.
\end{align*}
The $S_{\nu}$ are neither of finite dimension nor of finite codimension. However, they are semi-infinite in the sense that:
\[ \Gr^{\lambda}\cap S_{\nu} \mbox{ is of pure dimension $\langle \lambda + \nu, \rho^{\vee}\rangle$}, \]
where $2\rho^{\vee}$ is the sum of positive roots.

The varieties $\Gr^{\lambda}\cap S_{\nu}$ may be defined over finite fields. Let
\[ L_{\lambda\nu}(q) = |\{\mbox{$\FF_q$-rational points of $\Gr^{\lambda}\cap S_{\nu}$}\}|, \]
where $\FF_q$ denotes the finite field with $q$ elements.

We now work with symmetric functions in the group algebra (over $\ZZ[q^{-1}]$) of $X$. Here `symmetric' refers to invariant under the evident action of the Weyl group on $X$. 
The canonical reference for symmetric functions is \cite{Mac2}.

For $\lambda \in X_+$, let
\begin{align*}
m_{\lambda} &= \mbox{the monomial symmetric function corresponding to $\lambda$}, \\
P_{\lambda} &= \mbox{the Hall-Littlewood function corresponding to $\lambda$}.
\end{align*}
It is known that
\[ P_{\lambda} = \sum_{\nu\in X_+} q^{-\langle \lambda+\nu, \rho^{\vee}\rangle}L_{\lambda\nu}m_{\nu}. \]

\subsection{Soft Hodge theory}
Let $Y$ be a complex variety. The compactly supported (rational) cohomology $H_c^*(Y)$ admits a functorial (mixed) Hodge structure \cite{De}. Define
\[ h_c^{i,j;k}(Y) = \dim \, \gr^i_F\, \gr^W_{i+j}\, H_c^k(Y), \]
where $\gr_F$ (resp. $\gr^W$) denotes associated graded of the Hodge (resp. weight) filtration.
The compactly supported mixed Hodge polynomial and Hodge-Euler characteristic of $Y$ are defined by
\begin{align*}
h_c(Y; x,y,t) &= \sum h_c^{i,j;k}x^iy^jt^k,\\
E(Y; x,y) &= h_c(Y; x,y, -1),
\end{align*}
respectively. The Hodge-Euler characteristic is motivic: if $Y = \bigsqcup_{i=1}^n Z_i$, with each $Z_i$ a locally closed subvariety, then
\[ E(Y;x,y) = \sum_{i=1}^n E(Z_i; x,y).\]
Similarly, the K\"unneth formula yields:
\[ E(Y\times Z; x,y) = E(Y;x,y)E(Z;x,y).\]
If each $H^k_c(Y)$ is Hodge-Tate, i.e., $h_c^{i,j;k}(Y)=0$ for $i\neq j$, then $h_c(Y;x,y,t)$ is a polynomial in $xy$ and $t$, and we set
\begin{align*}
h_c(Y; q,t) &= h_c(Y; q^{\frac{1}{2}}, q^{\frac{1}{2}}, t) \\
E(Y; q) &= E(Y; q^{\frac{1}{2}}, q^{\frac{1}{2}}).
\end{align*}

\begin{scholium}
Hall-Littlewood polynomials encode the Hodge-Euler characteristics of the $\Gr^{\lambda}\cap S_{\nu}$. More precisely, for $\lambda\in X_+$,
\[ P_{\lambda} = \sum_{\nu\in X_+} q^{-\langle \lambda + \nu, \rho^{\vee}\rangle} E(\Gr^{\lambda}\cap S_{\nu}; q)m_{\nu}.\]
\end{scholium}

\begin{proof}
By \cite[Lemma 4]{GL}, each intersection $\Gr^{\lambda}\cap S_{\nu}$
admits a decomposition into locally closed subvarieties isomorphic to a product of some linear affine space ($\AA^s$) and some algebraic torus ($\GG_m^t$). In particular, $H^k_c(\Gr^{\lambda}\cap S_{\nu})$ is Hodge-Tate. Further, the Hodge-Euler characteristic is motivic, and
\begin{align*}
E(\AA^s\times \GG_m^t; q) &= q^s(q-1)^t \\
&= |\{\mbox{$\FF_q$-rational points of $\AA^s \times \GG_m^t$}\}|. \qedhere
\end{align*}
\end{proof}
I would be very grateful if someone could explain to me how to compute the polynomials $h_c(\Gr^{\lambda}\cap S_{\nu}; q,t)$. Apart from intrinsic interest, the problem of computing $H_c^*(\Gr^{\lambda}\cap S_{\nu})$ is similar in flavor to understanding certain intersections in the (finite) flag variety (see \cite{V}).

\end{document}